\UseRawInputEncoding
\documentclass{amsart}
\usepackage{amssymb}

\newtheorem{theorem}{Theorem}[section]

\theoremstyle{definition}

\newtheorem{corollary}[theorem]{Corollary}
\newtheorem{proposition}[theorem]{Proposition}
\newtheorem{remark}[theorem]{Remark}

\numberwithin{equation}{section}
\newcommand{\beqa}{\begin{eqnarray*}}
\newcommand{\eeqa}{\end{eqnarray*}}
\newcommand{\beqn}{\begin{eqnarray}}
\newcommand{\eeqn}{\end{eqnarray}}

\newcommand{\iy}{\infty}

\newcommand{\ov}{\overline}

\newcommand{\ci}{\subseteq}
\newcommand{\pf}{\noindent {\bf Proof :} }

\renewcommand{\a}{\alpha}
\renewcommand{\b}{\beta}

\newcommand{\e}{\varepsilon}

\newcommand{\m}{\mu}

\newcounter{cnt1}
\newcounter{cnt2}
\newcounter{cnt3}
\newcommand{\blr}{\begin{list}{$($\roman{cnt1}$)$}
        {\usecounter{cnt1} \setlength{\topsep}{0pt}
                \setlength{\itemsep}{0pt}}}
\newcommand{\bla}{\begin{list}{$($\alph{cnt2}$)$}
        {\usecounter{cnt2} \setlength{\topsep}{0pt}
                \setlength{\itemsep}{0pt}}}
\newcommand{\bln}{\begin{list}{$($\arabic{cnt3}$)$}
        {\usecounter{cnt3} \setlength{\topsep}{0pt}
                \setlength{\itemsep}{0pt}}}
\newcommand{\el}{\end{list}}
\newtheorem{thm}{Theorem}

\newtheorem{cor}[thm]{Corollary}

\newtheorem{Def}[thm]{Definition}
\newtheorem{prop}[thm]{Proposition}
\newtheorem{rem}[thm]{Remark}
\newcommand{\Rem}{\begin{rem} \rm}
\newcommand{\bdfn}{\begin{Def} \rm}
\newcommand{\edfn}{\end{Def}}

\title{On Ball dentable property  in Banach Spaces}
\author{Sudeshna Basu}
\address{Department of Mathematics \\George Washington University\\ Washington DC 20052\\ USA} 
\email{sbasu@gwu.edu, sudeshnamelody@gmail.com}
\date{}
\begin{document}
\maketitle
\begin{abstract}
In this work, we introduce the notion of  Ball dentable property in Banach spaces. We study certain stability results for the $w^*$-Ball dentable property
leading to a discussion on Ball dentability  in the context of ideals of Banach spaces.  We  prove that the $w^*$-Ball-dentable property
can be lifted from an $M$-ideal to the whole Banach Space. We also prove similar results for strict ideals  of a Banach space.
 We note that the space  $C(K,X)^*$ has $w*$-Ball dentable property when $K$ is dispersed and $X^\ast$ has the $w^*$-Ball dentable property.
\end{abstract}
\section{Introduction}
Let $X$ be a {\it real} Banach space and $X^*$ its dual. We will denote by $B_X$, $S_X$ and $B_X(x, r)$ the closed unit ball, the unit sphere and the closed ball of radius $r >0$ and center
$x$. We refer to the monograph \cite{B}
 for notions of convexity theory that we will be using here.
\bdfn
\blr
\item We say $A \ci B_{X^*}$ is a norming set for $X$ if $\|x\| =
\sup\{|x^*(x)| : x^* \in A\}$, for all $x \in X$. A closed subspace $F
\ci X^*$ is a norming subspace if $B_F$ is a norming set for $X$.
\item Let $f \in X^*$, $\a > 0$ and $ C \subseteq X$.
Then the set $S(C, f, \a) = \{x \in C : f(x) > \mbox{sup}~ f(C) - \a \}$ is called the open slice determined by $f$ and
$\a.$  We assume without loss of generality that $\|f\| = 1$. One can analogously define $w^*$ slices in $X^*$.
\item
A point $x \neq 0$ in a convex set $K \ci X$ is called a denting point 
point of $K$, if for every $\e > 0$, there exist slices
$S$ of $K$, such that $x \in S$ and $diam(S) <\e.$ One can analogously define $w^*$-denting point in $X^*$.
\el
\edfn
The following result will be useful in our discussion.
\begin{proposition}\cite{LLT}\label{llt}
Suppose $x \in B_X$. Then $x$ is denting point if and only if $x$ is a PC(point of continuity,i.e. points  for which the identity mapping on
the unit ball, from weak topology to norm topology is continuous.)  and an extreme point of $B_X.$
\end{proposition}
\begin{remark}
The above result is also true for $w^*$-denting points.
\end{remark}

\bdfn A Banach Space is said to have Ball-Dentable property (BDP)
if the unit ball has slices of arbitrarily small diameter. Analogously we
can define $w^*$-Ball dentability in a dual space. \edfn 
\begin{rem}
 \blr
\item It is clear that if $B_X$ has a denting points point, then it has BDP.
\item Clearly Radon Nikodym Property (RNP) implies BDP.
( see \cite{B}.)
\el
\end{rem}

In this short note, we study certain stability results for $w^*-$BDP
leading to a discussion on BDP in the context of ideals of Banach
spaces, 
see \cite{GKS} and \cite{R}. We use various techniques from
the geometric theory of Banach spaces to achieve this. The spaces
that we will be considering have been well studied in the literature.
 A large class of function spaces like the Bloch spaces, Lorentz
and  Orlicz spaces, spaces of vector-valued functions  and spaces of
compact operators are examples of the spaces we will be considering, for details, see \cite{HWW}.
We provide some descriptions of $w^*$-denting  points in Banach spaces in different contexts.
We need the following definition.
 \bdfn Let $X$ be a Banach space.
\blr
\item A linear projection $P$ on $X$ is called an  M-{\it projection} if \\
$ \|x \| = max \{ \|Px\|, \|x - Px\| \},$ for all $ x \in X$.\\
A linear projection $P$ on $X$ is called an  L-{\it projection} if \\
$\|x\| =   \|Px\| + \|x - Px\| $ for all $ x \in X$.
A linear projection $P$ on $X$ is called an  $L^{p}$-{\it projection}($ 1< p < \infty $) if \\
$\|x\|^{p} =   \|Px\|^{p} + \|x - Px\|^{p} $ for all $ x \in X$.
\item A subspace $M \subseteq X$ is called an $M$-summand if it is the range of an $M$-projection. A subspace $M \subseteq X$ is called an $L$-summand if it is the range of an $L$-projection. A subspace $M \subseteq X$ is called an $L^{p}$-summand if it is the range of an $L^{p}$-projection.
\item A closed subspace $M \subseteq X$ is called an $M$-ideal if $M^{\perp}$ is the kernel of an $L$-projection in $X^*$.
\el
\edfn
We recall from Chapter I of \cite{HWW} that when $M \subset X$ is an
$M$-ideal, elements of $M^\ast$ have unique norm-preserving
extension to $X^\ast$ and one has the identification, $X^\ast =
M^\ast \oplus_1 M^\perp$.
 Several examples from among
function spaces and spaces of operators that satisfy these geometric
properties can be found in the monograph \cite{HWW}, see also \cite{E}.
 First, we prove for an $M$- ideal $M \subset X$, any $w^\ast$-denting  point of $B_{M^\ast}$,is also a $w^\ast$-denting  point of $B_{X^\ast}$ . We prove similar  for  a strict ideal $Y \subset X$(see Section 2 for the definition) i.e., we prove that a $w^*$-denting  point of a strict ideal continues to be so in the bigger space. We also prove corresponding results for the Ball dentable property. The techniques used in the proofs are adapted from \cite{BR}.

\section{Stability Results}
We will use the standard notation of $\oplus_1$,$\oplus_p$,and  $\oplus_{\iy}$ to denote the
$\ell^1 \ell^p,$ and $\ell^\infty$-direct sum of two or more Banach spaces.
Let $M \subseteq X$ be an $M$-ideal. It follows from the results in Chapter I in \cite{HWW}
that any $x^\ast \in X^\ast$, if $\| m^\ast\|= \|x^\ast|M\| = \|x^\ast\|$, then $x^\ast$ is the unique norm preserving extension of $m^\ast$. For notational convenience we denote both the functionals by $m^\ast$. Clearly any $M$-ideal is also an ideal. 

\begin {proposition} 
Let $Z = X\oplus_{1}Y.$ Let $x_0$ be a denting point of $B_X.$ Then $x_0$ is a denting point of $B_Z.$
\end{proposition}
\pf
Let $\{z_n\}$ be a sequence in $B_Z$ such that $z_n \longrightarrow x_0,$ weakly.If P denotes the projection mapping to $X$, it follows that $P(z_n)\longrightarrow x_0.$ Since $x_0$ is a denting poing of $B_X,$ it follows from Proposition \ref{llt}, that $x_0$ is a PC. So  $P(z_n)\longrightarrow x_0$ in norm.
Since $x_0$ is denting, it is an extreme point as well ,so $\|x_0\|=1,$ and it follows that $\underline{lim}P(x_n)=1$ and $\underline{lim}(x_n)=1.$
This implies $\lim_{n}\| z_n -P(z_n)\| = 0.$
Therefore $\|z_n- x_0\| \leq \|z_n- P(z_n)\| +\|z_n- x_0\| \longrightarrow 0.$
Thus $x_0$ is a point of continuity of $B_Z$
Also, since $x_0$ is an extreme point of $B_X$ ,it is a an extreme point of $B_Z.$
Again using Proposition \ref{llt}we have, $x_0$ is a denting point of $B_Z.$\qed
\Rem
Similar conclusion follow for $w^*$ -denting points also.
\end{rem}
The following corollary is immediate.
\begin{corollary}
Suppose $M\subseteq X$ is an M-ideal. If $m^* \in B_{M^*}$ is a $w^*$-denting point of $B_{M^*},$ then $m^* \in B_{X^*}$ is a $w^*$-denting point of $B_{X^*}.$
\end{corollary}

In the case of an $M$-ideal $M \subset X$, the following is true.
\begin{prop} \label{M-idealBDP}
Let $M \subseteq X$ be an $M$-ideal, then if $M^*$ has the $w^*$-BDP then $X^*$ has the  $w^*$-BDP.
\end{prop}
\pf
 Suppose $M^*$ has the $w^*$-BDP, then for any for any
$\e>0$ there exists slices $S_{M}$ and  $S_{M}= \{ m^* \in B_{M^*} / m^*(m_0)> 1 - \a \}$
and $(dia S_{M}) < \e.$ Since $M$ is an $M$-
ideal, for any $x^*\in X^*$ we have the unique decomposition ,$x^*=
m^*+ m^{\perp}$, where $m^* \in M^\ast$ and $m^\perp \in M^\perp$.
Suppose we have $0<\m<\a.$ Then 
\beqa
S_{X} &= &\{ x^* \in B_{X^*} / x^*(m_0) > 1 - \mu\}\\
&= & \{x^* \in B_{X^*}/ m^*(m_0) +m^{\perp}(m_0) > 1 - \mu\}\\
&  \subseteq &S_M \times \mu B_{M^\perp}\\
 \eeqa 
 Choose
$\b = min ( \m, \e).$ Then 
\beqa
S'_{X}& =& \{ x^* \in B_{X^*} / x^*(m_0) > 1 - \b\} \subseteq S_{X} \times \b B_{M^\perp}\\
\eeqa
Thus $dia( S'_{X} )\leq  dia(S_{M}) + 2\e < \e+2\e = 3\e.$\\
Also, since $ \|m_0\|=1 ,$ there exists $m^*_{0} \in B_{M^*}$ such
that $m^*_{0} (m_0) > 1 - \b.$ Hence $m^*_{0}\in S'_{X}.$
 \qed

Arguing similarly it follows that 
 \begin{cor}\label{lp-BSCSP}
Suppose $X = \oplus_{p}X_i, 1 \leq p \leq \infty $. If
 $X_{i}^*$ has the $w^*$-BDP for some $i$, then $X^*$ has the $w^*$-BDP.
\end{cor}
 The above arguments extend easily to vector-valued
continuous functions. We recall that for a compact Hausdorff space
$K$, $C(K,X)$ denotes the space of continuous $X$-valued functions
on $K$, equipped with the supremum norm. We recall from \cite{La} that dispersed compact Hausdorff spaces
have isolated points.
\begin{cor}
Suppose $K$ is a compact Hausdorff space with an isolated point. If $X^*$ has the $w^*$-BDP, then $C(K,X)^*$ has
the $w^*$-BDP .
\end{cor}
\pf  Suppose $X^*$ has the $w^*$-BDP. For
an isolated point $k_0 \in K$, the map $F \rightarrow \chi_{k_0}F$
is a $M$-projection in $C(K,X)$ whose range is isometric to $X$. Hence we see that
$C(K,X)^*$ has the $w^*$-BDP. \qed

 \bdfn We recall that an ideal $Y$ is said to
be a strict ideal if for a projection $P: X^\ast \rightarrow X^\ast$
with $\|P\| = 1$, $ker(P) = Y^\bot$, one also has  $B_{P(X^\ast)}$ is $w^*$-dense in $B_{X^*}$ or
in other words $B_{P(X^*)}$ is a norming set for $X.$ \edfn
 In the case of an
ideal also one has that $Y^\ast$ embeds (though there may not be
uniqueness of norm-preserving extensions) as $P(X^\ast)$. Thus we
continue to write $X^* = Y^* \oplus Y^\perp$. In what follows we use
a result from \cite{R1}, that identifies strict ideals as those for
which $Y \subset X \subset Y^{\ast\ast}$ under the canonical
embedding of $Y$ in $Y^{\ast\ast}$.
\begin{prop}
Suppose $Y$ is a strict ideal of $X$. If $y^*\in B_{Y^*}$ is a $w^*$-denting point of $B_{Y^*},$ then $y^*$ is a $w^*$-denting point of $B_{X^*}.$
\end{prop}
\pf Since $y^*\in B_{Y^*}$ is a $w^*$-denting point of $B_{Y^*},$ for
any $\e>0$ there exists $w^*$ slices $S$ and  $S= \{ y^* \in B_{Y^*} / y^*(y_0)> 1 -
\a \}$ and $dia( S_) < \e.$ Since $Y$ is a
strict ideal in $X$ , we have $B_{X^*}= \ov{B_{Y^*}}^{w^*}$, hence
we have the following \beqa
 S' & = & \{ x^* \in B_{X^*} / x^*(x_0)> 1 - \a \}\\
&=&\{ x^* \in \ov{B_{Y^*}}^{w^*} / x^*(x_0)> 1 - \a \}\\
\Longrightarrow  dia( S') &= &dia(S) < \e \\
 \eeqa
 Hence $y^*$ is a $w^*$-denting  point of $B_{X^*}.$
\qed

Arguing similarly it follows that:
\begin {prop}
Suppose $Y$ is a strict ideal of $X$.
If $Y^*$ has $w^*$-BDP then $X^*$ has $w^*$-BDP.
\end{prop}
\begin{rem}
A prime example of a strict ideal is a Banach space $X$ under its canonical embedding in $X^{\ast\ast}$.
It is known that any $w^*$-denting point of $B_{X^{\ast\ast}}$ is a point of $X$. 
\end{rem}


\begin{thebibliography}{99}
\small
\bibitem [A]{A} M.\ D.\ Acosta, A.\ Kaminska and M. \ Mastylo, {\it
The Daugavet property in rearrangement invariant spaces}, Trans.
Amer. Math. Soc., {\bf 367} (2015) 4061--4078.

\bibitem[B]{B} R.\ D.\ Bourgin, {\it Geometric aspects of convex sets with the Radon-Nikodým property}, Lecture Notes in Mathematics, 993. Springer-Verlag, Berlin, 1983. xii+ 474 pp.
\bibitem[BR]{BR} S.\ Basu and T.\ S.\ S.\ R.\ K.\ Rao { \it On Small Combination of slices in Banach Spaces}, Extracta Mathematica, {\bf 31}( 2016) 1--10.  

\bibitem[GKS] {GKS} G.\ Godefroy, N.\ J.\  Kalton and P.\ D.\ Saphar,  {\it Unconditional ideals in Banach spaces}, Studia Math. {\bf 104}  (1993) 13-59.
\bibitem[LLT]{LLT} B.\ L.\ Lin , P.\ K.\ Lin and S.\ Troyanski, {\it Charecterisation of denting points}, Proc. Amer.
Math. Soc. { \bf 102} (1988) 526--528.
\bibitem [HWW]{HWW} P.\ Harmand, D.\ Werner and W.\ Werner, {\it $M$-ideals in Banach spaces
and Banach algebras}, Lecture Notes in Mathematics, 1547.
Springer-Verlag, Berlin, 1993. viii+387 pp.
\bibitem[La]{La} H. E. Lacey, {\it The isometric theory of classical Banach
spaces}, Die Grundlehren der mathematischen Wissenschaften, Band
208. Springer-Verlag, New York-Heidelberg, 1974. x+270 pp.
\bibitem [E]{E} E.\ Oja, {\it On $M$-ideals of compact operators in
Lorentz sequence spaces}, Journal of Mathematical Analysis and
Application, {\bf 259} (2001) 439--452.
\bibitem[R1]{R1}  T.\ S.\ S.\ R.\ K.\ Rao, {\it On ideals in Banach spaces}, Rocky Mountain J. Math.  {\bf 31} (2001) 595-609.
\bibitem[R]{R} H. \ P.\ Rosenthal, {\it On the structure of non-dentable closed bounded convex sets}, Adv. in Math.,  { \bf 70} (1988) 1--58.

\end{thebibliography}
\end{document}